# An Efficient Interval Uncertainty Optimization Approach Based on Quasi-sparse Response Surface

**Kefeng Wang[1], Pu Li[2], Yanfeng Zhang[3], Yunbao Huang[1]**
[1]Provincial Key Laboratory of Computer Integrated Manufacturing, Guangdong University of Technology, Guangzhou, Guangdong, China;
[2]School of Physics and Electrical Engineering, Shaoguang University, Shaoguan, Guangdong, China;
[3]School of Electronic and Information Engineering, Harbin Institute of Technology, Shenzhen, P.R. China

## ABSTRACT

The structure uncertainty optimization problem is usually treated as double-loop optimization process, which is computation-intensive. In this paper, an efficient interval uncertainty optimization approach based on Quasi-sparse response surface (QSRS) is proposed for structure uncertainty optimization. In which, 1) with $l_1$ norm and $l_2$ norm penalty method, a few appropriate basis functions are selected form a large number of basis functions to construct QSRS accurately and only a few sampling points is required, 2) as the orthogonal chebyshev polynomials is employed as QSRS basis functions, the local uncertainty can be evaluated by the combination of QSRS and interval arithmetic. Hence, the inner optimization process is eliminated. One mathematical problem and one engineering problem are used to validate the efficiency of the proposed approach. The results show that only 25% sampling points is needed than recently published approach.

## 1. INTRODUCTION

In structure design optimization, there are many uncertain parameters (such as material property) that may cause significant performance variation. Hence, an appropriate optimization method should consider these parameters. In practice, many uncertain parameters are uncertainty-but-bounds parameters and their bounds can be easily obtained. Thus, the convex model and interval model, in which the uncertainty-but-bounds parameters are treated as interval number, is proposed for structure uncertainty optimization. Interval models often involve a nest double-loop optimization process, it means that the simulation model in the outer loop optimization process should be evaluated iteratively by an inner optimization process [1]. It's costly, especially for complex simulation models (i.e. Finite element model).

In order to reduce the computational cost, the response surface is considered to approximate the simulation model. The matrix perturbation theory through a first order Taylor series expansion was applied to obtain a conservative dynamic response of interval functions [2]. The first-order Taylor series expansion was used to simplify the double loop optimization process [3]. Chen [4] used first order Taylor series expansion to calculate the upper bound of response for the robustness analysis of the interval vibration control system. However, because of low accuracy of the first-order Taylor series, the final optimization results obtained by those approaches may be in the infeasible region. Later, the interval optimization approach combining chebyshev response surface and interval arithmetic is proposed by Wu et al [5]. In this approach, the uncertainty caused by the uncertainty-but-bounds parameters is also treated as interval number. In the outer loop optimization process, Multi-Island genetic algorithm (MIGA) is employed to search global optimal points. In the inner loop optimization process, the chebyshev response surface is constructed in the uncertain range of design points updated in the outer loop optimization process, then, the bounds of uncertainty can be easily obtained by interval arithmetic. Hence, the inner loop optimization process is eliminated, which means the computation cost is greatly reduced. However, the sampling points required in chebyshev response surface is multiplied exponentially with the expansion order and the function dimension. For high-dimension engineering problems, the sampling points will be increased dramatically. Thus, the chebyshev collocation method (CCM) [6] is employed to reduce the sampling size, but the expansion order of chebyshev polynomials will be decreased, which means that the accuracy of chebyshev response surface may be decreased. What's more, the sampling size is still too large.

In this paper, an efficient interval uncertainty optimization approach using QSRS is proposed for structure uncertainty optimization (so-called QSRS approach). The QSRS [7], which can achieve a high accuracy with smaller sampling size than other response surfaces, is introduced for structure uncertainty optimization. In QSRS, a large number of orthogonal chebyshev polynomials are used as basis functions to enhance its expression ability. With $l_1$ norm and $l_2$ norm penalty method, only a few basis functions which are important to explain the relationship between variables and response values are selected, therefore, the coefficients associated with the selected basis functions can be estimated with a few sampling points. Additionally, the uniform design is employed as sampling method to sample the complete feature of source model. For the optimization process, MIGA is employed to search global optimum and update design point, then, the interval arithmetic combining with QSRS is used to directly calculate the bounds of objective in the uncertainty range of design point. Hence, the inner loop optimization process is eliminated. The significances of this approach include:

(1) The number of basis functions in QSRS will be increased to enhance its expression ability.
(2) With $l_1$ norm and $l_2$ norm penalty method, only a few basis functions are selected from a large number of basis functions to construct QSRS, thus its generalization ability and robustly are improved and less sampling points are required.
(3) The orthogonal chebyshev polynomials are selected as QSRS basis functions, the maximum value of object function can be evaluated by interval arithmetic, which means the inner optimization process is eliminated.

The rest of the paper is organized as follow: Section 2 introduces the uncertainty optimization problem. Section 3 describes QSRS approach. Then, we test the performance of QSRS approach and discuss the experiment results in section 4. Finally, we give the summary of some significant conclusions.

## 2. Uncertainty optimization problem

In some structure design problems, the design variables are also under uncertain. Hence, the uncertainty optimization problem is given in this section, in which the uncertain design variables and uncertainty parameters are both treated as interval numbers. Furthermore, the uncertainties of objective and constraints caused by these interval numbers can be evaluated.

The general structure design optimization problem with deterministic design variables and parameters is given as following:

$$\begin{cases} \min_{x} \quad f(x, y) \\ s.t. \quad g_i(x, y) \leq 0, i = 1,2,\cdots,n \\ \quad x^L \leq x \leq x^R \end{cases} \quad (1)$$

Where $x$ and $y$ are the vector of design variables parameters, respectively.

Any interval number can be defined as follows:
$$[x] = [\underline{x}, \bar{x}] = x_c + [-\Delta x, \Delta x] \quad (2)$$
where $\underline{x}$ and $\bar{x}$ are the bounds of $[x]$, $x_c = (\underline{x} + \bar{x})/2$ is the midpoint of $[x]$, $\Delta x = (\underline{x} - \bar{x})/2$ is the width of $[x]$.

When the uncertainty of design variables and parameters are considered, the problem (1) is transform as:

$$\begin{cases} \min_{x} \quad f([x], [y]) \\ s.t. \quad g_i([x], [y]) \leq 0, i = 1,2,\cdots,n \\ \quad x^L \leq x \leq x^R \end{cases} \quad (3)$$

where the bounds of uncertain parameters $[y]$ are usually pre-given, the width of uncertain design variables $[x]$ is also pre-given as $\xi$, which is defined as:
$$[x] = x_c + [-\xi, \xi] \quad (4)$$
Hence, the objective $f$ and constraints $g_i$ in problem (3) both are interval numbers and are written as $[f]$ and $[g_i]$, respectively. The problem of minimizing $[f]$ can be regarded as minimizing both the midpoint $f_c$ and width $\Delta f$, which is similar to the traditional robust design optimization [8]. Then, the problem (3) can be seen as multi-objective optimization problem. Thus, the weigh method is used to solve this problem, as the weight factors of midpoint $f_c$ and width $\Delta f$ both are set as 1, the objective of problem (3) is transformed to $f_{obj} = f_c + \Delta f$ which is the upper bounds of $[f]$.

For the constraints, in order to ensure the design points is within the feasible region, the upper bounds of are $[g_i]$ is employed to meet constraints:
$$\bar{g}_i([x], [y]) \leq 0 \quad (5)$$
Finally, the problem (3) is reformulate as:
$$\begin{cases} \min_{x_c} \max_{x \in [x], y \in [y]} f(x, y) \\ s.t. \max_{x \in [x], y \in [y]} g_i(x, y) \leq 0, i = 1,2,\cdots,n \\ \quad x^L \leq x \leq x^R \end{cases} \quad (6)$$

The problem (6) is usually solved by double-loop optimization process, however it's computation-consumption. For instance, if 1000 iterations of optimization process are required in the outer loop, 1000 iterations optimization process are required in inner loop, the total numbers of simulation evaluation are $1000 \times 1000 = 10^6$, which is costly. Hence, an efficient QSRS approach is proposed for the problem (6). In the next section, we will introduce the proposed QSRS approach in details.

## 3. QSRS approach

The QSRS approach is based on QSRS and interval arithmetic, in which the QSRS is based on the sparse representation theory. Hence, we will introduce sparse representation theory, SRS structure, elastic net and optimization process in turn.

### 3.1 Sparse representation

The sparse representation theory reveals that a function or a signal can be represented by linear combination of a few significant terms from a set of basis functions. It means only a few coefficients corresponding to the basis functions are large in magnitude, while most of them are zero or close to zero. Suppose the function $f \in \mathbb{R}^{m \times 1}$ has a sparse representation on the basis matrix $\mathbf{\Phi} \in \mathbb{R}^{m \times n}$ (see section 3.2), we will have:
$$f \approx \mathbf{\Phi}\boldsymbol{\beta} \quad (7)$$
where $\boldsymbol{\beta} \in \mathbb{R}^{n \times 1}$ is the vector of coefficients whose elements are almost zero.

In order to better demonstrate this theory, consider an 2-dimension function: $g(x) = (x_1 + 2x_2 - 7)^2 + (2x_1 + x_2 - 5)^2, -10 \leq x_1, x_2 \leq 10$. Suppose $g(x)$ can be represented by $\mathbf{\Phi}$, that is, $\forall$ a small enough number $\varepsilon \geq 0, \exists c_i \in \mathbb{R}$ such that:
$$\max_{-10 \leq x, y \leq 10} \left| g(x) - \sum_{i=1}^{n} c_i \varphi_i(x) \right| \leq \varepsilon \quad (8)$$

Where $\varphi_i(x)$ is the basis functions, which is Chebyshev polynomials and their number is set as $n = 30$. Their corresponding coefficients can be calculated as follows:
$$c_i = \frac{2}{\pi} \int_0^\pi \frac{g(x)\varphi_i(x)}{\sqrt{1-x^2}} dx = \frac{2}{\pi} \int_0^\pi g(\cos\theta)\varphi_i(\theta) d\theta \quad (9)$$

Note that the function domain should be mapped to domain $[0, \pi]$ (see section 3.2). The solution of coefficients is shown in figure 1, it can be seen that only 6 coefficients are large in magnitude and the rest are close to zero. The figure of actual function and reconstructed function is drawn by sampling $50 \times 50$ points in function domain which is shown in figure 2 and figure 3, respectively. We can observe that the function $g(x)$ is well appropriated by 6 Chebyshev polynomials. The sparseness of $g(x)$ in Chebyshev polynomials is 6.

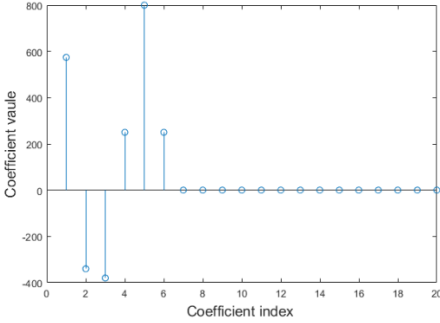

**FIGURE 1:** COEFFICIENTS CALCULATED BY Eq (9)

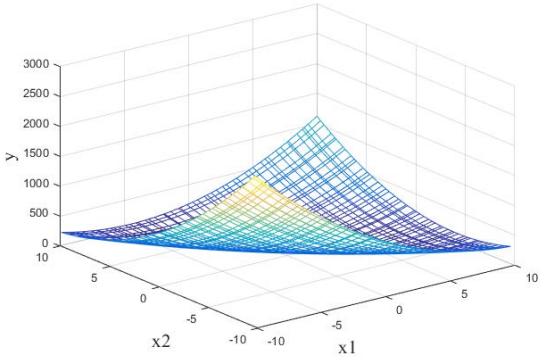

**FIGURE 2:** ACTUAL FUNCTION

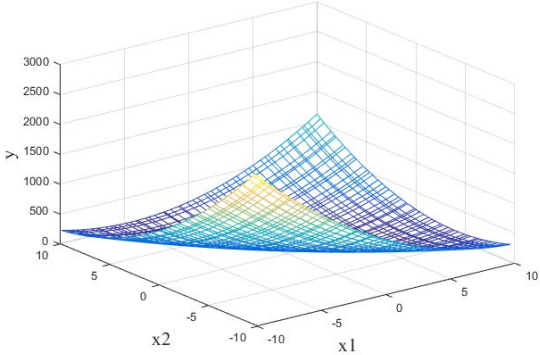

**FIGURE 3:** RECONSTRUCTED FUNCTION

### 3.2 QSRS structure

The QSRS is based on polynomial response surface which can be written as follows:

$$f(x) \approx \hat{f}(x) = \sum_{i=1}^{n} \varphi_i(x) \beta_i \quad (9)$$

Where $x = [x^{(1)}, x^{(2)}, \cdots x^{(p)}]$ is a sampling point which is obtained by uniform design method, $p$ is the number of variables. $\beta_i$ is the corresponding coefficient of basis function.

$\varphi_i(x), i = 1, 2, \cdots n$ are the basis functions (so-called atoms), $n$ is the number of basis functions. In this paper, the Chebyshev polynomials are selected as basis functions, as the wrapping effect of interval arithmetic with Chebyshev polynomials can be well monitored. Suppose a variable $x \in [-1,1]$, the chebyshev polynomials of $x$ is defined as $\varphi_i(x) = \cos(i\theta), \theta = \arccos(x), \theta \in [0, \pi]$. For $m$-dimensional cases, the basis functions are the tensor product of basis functions in one-dimension.

The Eq (9) can also be written as matrix form:

$$\hat{f} = \Phi \beta \quad (10)$$

where $\Phi$ is design matrix (so-called dictionary) which is consist of atoms. Given a set of sampling points $x_i = [x_1, x_2, \cdots x_m]$, $\Phi$ can be defined as:

$$\Phi = \begin{bmatrix} \varphi_1(x_1) & \cdots & \varphi_n(x_1) \\ \vdots & \ddots & \vdots \\ \varphi_1(x_m) & \cdots & \varphi_n(x_m) \end{bmatrix} \quad (11)$$

$\beta$ is the vector of corresponding coefficients. Note that the sampling points $x = [x_1, x_2, \cdots x_m]$ should be mapped to $[0, \pi]$ and the number of atoms is set as $p = 3m$.

Motivated by sparse represent theory, we will seek sparse representation of $f$ in dictionary $\Phi$, which means that the process of constructing QSRS is to solve the problem as follows:

$$\|f - \Phi\beta\|_2^2 \leq \varepsilon \ s.t. \|\beta\|_0 \quad (12)$$

where $s = \|\beta\|_0$ is the sparseness of $f$ in dictionary $\Phi$.

Since the $l_0$ norm problem is a non-convex NP-hard problem, we relax it to $l_1$ norm problem as follows:

$$\|f - \Phi\beta\|_2^2 \leq \varepsilon \ s.t. \|\beta\|_1 \quad (13)$$

the $l_1$ norm problem can also be written as Lagrange multiplier form as follows:

$$\hat{\beta}(\lambda) = \underset{\beta}{argmin} \ \|f - \Phi\beta\|_2^2 + \lambda\|\beta\|_1 \quad (14)$$

It's the famous LASSO (least absolute shrinkage and selection operator) model. The LASSO regression can reconstruct the source model accurately in the case of the sparseness of source model is smaller than the number of sampling points. However, the number of atoms selected by LASSO regression cannot exceed the number of sampling points. In practice, the source model may be too complex so that its sparseness in the selected basis functions is not small enough. In this case, the response surface with LASSO regression will be unstable. In addition, we can increase the number of sampling points, but it's costly. Hence, the elastic net regression which will increase the number of selected atoms is introduced to construct the QSRS.

### 3.3 Elastic net

Suppose a function combining with $l_1$ norm penalty, $l_2$ norm penalty and cost function as follows:

$$\mathcal{L}(\lambda_1, \lambda_2, \beta) = \|f - \Phi\beta\|^2 + \lambda_1\|\beta\|_1 + \lambda_2\|\beta\|_2^2 \quad (12)$$

where

$$\|\beta\|_1 = \sum_{j=1}^{n} |\beta_j|$$

$$\|\boldsymbol{\beta}\|_2^2 = \sum_{j=1}^n \beta_j^2$$

The Elastic net regression aims to find the minimal solution $\hat{\boldsymbol{\beta}}$ of Eq (12):

$$\hat{\boldsymbol{\beta}} = argmin\{\mathcal{L}(\lambda_1, \lambda_2, \boldsymbol{\beta})\} \quad (13)$$

where the $l_2$ norm penalty denotes LASSO regression, $l_2$ norm penalty denotes ridge regression. To better illustrate their relationship with elastic net, let $\alpha = \lambda_2/(\lambda_1 + \lambda_2)$, the Eq (12) can be written as follows:

$$\mathcal{L}(\lambda, \alpha, \boldsymbol{\beta}) = \|\boldsymbol{f} - \boldsymbol{\Phi}\boldsymbol{\beta}\|^2 + \lambda((1-\alpha)\|\boldsymbol{\beta}\|_1 + \alpha\|\boldsymbol{\beta}\|_2^2) \quad (14)$$

The parameter α determines the mixing ratio of the ridge regression and the LASSO regression. When $\alpha = 1$, the elastic net becomes ridge regression, highly correlated atoms tend to be selected at the same time. When $\alpha = 0$, it becomes LASSO regression, a sparse solution will be given in this case.

The Eq (12) has two parameters, which can be solved by referring to the following theorems.

**Theorem 1.** Given data set $(\boldsymbol{y}, \boldsymbol{X})$ and $(\lambda_1, \lambda_2)$, the elastic net estimation is:

$$\boldsymbol{\beta}(elastic\ net) = \underset{\boldsymbol{\beta}}{argmin}\ \boldsymbol{\beta}^T \left(\frac{\boldsymbol{X}^T\boldsymbol{X} + \lambda_2 \boldsymbol{I}}{1 + \lambda_2}\right)\boldsymbol{\beta} - 2\boldsymbol{y}^T\boldsymbol{X}\boldsymbol{\beta} + \lambda_1|\boldsymbol{\beta}|_1 \quad (15)$$

The LASSO solution is:

$$\boldsymbol{\beta}(lasso) = \underset{\boldsymbol{\beta}}{argmin}\ \boldsymbol{\beta}^T(\boldsymbol{X}^T\boldsymbol{X})\boldsymbol{\beta} - 2\boldsymbol{y}^T\boldsymbol{X}\boldsymbol{\beta} + \lambda_1|\boldsymbol{\beta}|_1 \quad (16)$$

It can be seen that the elastic net shrinks the correlation matrix $\boldsymbol{X}^T\boldsymbol{X}$ towards the identity matrix.

**Lemma 2** [9]. Given data set $(\boldsymbol{y}, \boldsymbol{X})$ and $(\lambda_1, \lambda_2)$. Define a new set of data sets $(\boldsymbol{y}^*, \boldsymbol{X}^*)$:

$$\boldsymbol{X}^* = (1+\lambda_2)^{-1/2}\begin{pmatrix}\boldsymbol{X}\\ \sqrt{\lambda_2}\boldsymbol{I}\end{pmatrix}, \boldsymbol{y}^* = \begin{pmatrix}\boldsymbol{y}\\\boldsymbol{0}\end{pmatrix} \quad (17)$$

Let $\gamma = \lambda_1/\sqrt{(1+\lambda_2)}$, $\boldsymbol{\beta}^* = \sqrt{(1+\lambda_2)}\boldsymbol{\beta}$, the LASSO problem is equaled with

$$\mathcal{L}(\gamma, \boldsymbol{\beta}) = \mathcal{L}(\gamma, \boldsymbol{\beta}^*) = |\boldsymbol{y}^* - \boldsymbol{X}^*\boldsymbol{\beta}^*|^2 + \gamma|\boldsymbol{\beta}^*|_1 \quad (18)$$

let $\boldsymbol{\beta}^* = \underset{\boldsymbol{\beta}^*}{argmin}\{\mathcal{L}(\gamma, \boldsymbol{\beta})\}$, then $\boldsymbol{\beta}(elastic\ net) = \sqrt{1+\lambda_2}\boldsymbol{\beta}^*$.

The proof of theorem and lemma are a simple derivation, not listed in the text. It can be seen from the lemma that the elastic net problem can be transformed to LASSO problem with the fixed parameter $\lambda_2$, and the LASSO solution can be obtained with the least angle regression (LAR). The MATLAB toolbox provided by reference [10] is employed to directly find the elastic net solution with fixed $\lambda_2$.

For the selection of the parameter $\lambda_2$, the following principles should be followed: 1) If the source model has sparse representation on the basis function, $\lambda_2$ should be 0, the elastic net problem is transformed to LASSO problem, 2) when the correlation of dictionary is not low, $\lambda_2 > 0$, and its value will be slowly increased as correlation of dictionary become higher. The value of $\lambda_2$ is generally smaller.

In this paper, the cross-validation method is employed to determine the value of the parameter $\lambda_2$. The specific process is shown in Table 1.

Table 1. The process of constructing QSRS

| | |
|---|---|
| Step 1 | Obtain sampling points $\boldsymbol{x}$ by uniform design. |
| Step 2 | Calculate response value $\boldsymbol{y}$ of $\boldsymbol{x}$ in the source model. |
| Step 3 | Divide the cross-validation sets. Divide the sampling points $\boldsymbol{x}$ and the response value $\boldsymbol{y}$ into k sets. k usually takes 10, or other values divided by the number of sampling points. |
| Step 4 | Fix $\lambda_2$ to 0, 0.0001, 0.001, 0.01, 0.1, 1, and 10, respectively, perform steps 5 to 7. |
| Step 5 | Take k-1 sets as the training set, and the remaining one is the prediction set. Use the training set to construct the QSRS by LAR-EN solver and evaluate its prediction accuracy on the prediction set. |
| Step 6 | Repeat step 5 k times to ensure that each set is used for prediction |
| Step 7 | Take the average prediction accuracy of k times experiments for every fixed $\lambda_2$. |
| Step 8 | Output the minimum $\lambda_2$. |
| Step 9 | Construct the QSRS by all sampling points for the fixed $\lambda_2$. |

### 3.4 Interval optimization process

Traditional double loop optimization includes: (1) the outer loop searches the global optimum in the design domain and updates design point, (2) the inner loop finds the maximum value of objective in the uncertain range of design point. In general, some intelligent algorithms are employed in outer loop optimization process and inner loop optimization process, however, it's costly. For QSRS approach, as the interval arithmetic has high computational efficiency, the QSRS combined with interval arithmetic is used to replace the costly inner loop optimization process.

The interval arithmetic refers to some operation rules between two interval numbers. The general operation rules between them is given as follows:

$$\begin{cases}[x] + [y] = [\underline{x} + \underline{y}, \overline{x} + \overline{y}]\\ [x] - [y] = [\underline{x} - \overline{y}, \overline{x} - \underline{y}]\end{cases} \quad (19)$$

As $\cos(\boldsymbol{\theta}) \in [-1,1]$ which is the basis function of QSRS, the bounds of $\hat{y}(x)$ can be calculated as follows:

$$[\hat{y}] = [\beta_0 - \sum_{i=1}^n |\beta_i|, \beta_0 + \sum_{i=1}^n |\beta_i|] \quad (20)$$

To better introduce the interval arithmetic used in this paper. Consider the function $y = x^3 + 2x, x \in [-1,1]$, the QSRS is used to approximate this function and its expression is given as follows:

$$\hat{y}(x) = 2.75\cos(\theta) + 0.25\cos(3\theta), \theta = \arccos(x) \quad (20)$$

Based on the Eq (20), the interval range of $\hat{y}$ is given as follows:

$$[\hat{y}] = 2.75 \times [-1,1] + 0.25 \times [-1,1] = [-3,3] \quad (21)$$

The actual range of $y$ in domain is also $[-3,3]$. Hence, the QSRS combined with interval arithmetic can be used to directly calculate the maximum value of objective in the uncertain range of design point. In this case, the inner loop optimization process is eliminated and the efficiency is improved.

To search the global optimum, MIGA is employed for the outer loop optimization process, and the settings of sub-population size in MIGA may affect its performance. Thus, the sub-population size will be adjusted for different problems.

The specific process of QSRS approach is shown in figure 4. Such method contains four main steps: (1) Initiation. $x_c$ denotes the mid value of design point, $\xi$ is the width of design point, $x^L, x^R$ is the domain of design variables. (2) Construct QSRS. First, we will generate some sampling points by uniform design. Then, calculate the response values of sampling points. Finally, calculate the coefficients of basis functions and construct QSRS. (3) Calculate the objective $f_{max}(x_c, y)$ and constraints $g_{max}(x_c, y)$ by interval arithmetic. (4) Update the design point by MIGA, if the result satisfies converge conditions, the optimization process will be end, otherwise it will go to step (2).

## 4. Experiments and analysis

In the section, a mathematical example and an engineering example is used to demonstrate the efficiency of QSRS approach.

### 4.1 Mathematical example

The "Three-Hump" function which is a test function for optimization is employed as a mathematical example, which is given by:
$f(x_1, x_2) = 2x_1^2 - 1.05x_1^4 + \frac{1}{6}x_1^2 + x_1 x_2 + x_2^2$
$-5 \leq x_1 \leq 5, -5 \leq x_2 \leq 5$

We assume that both the two design variables are interval design variables, and the widths of them are 0.1. We will have the uncertain optimization problem as follows:
$\min_{x_{1c}, x_{2c}} \max_{x_1 \in [x_1], x_2 \in [x_2]} f(x_1, x_2)$
s.t. $[x_1] = x_{1c} + [-\xi_1, \xi_1], [x_2] = x_{2c} + [-\xi_2, \xi_2]$
$-5 + \xi_1 \leq x_{1c} \leq 5 - \xi_1, -5 + \xi_2 \leq x_{2c} \leq 5 - \xi_2$

The interval optimization approach and QSRS approach are used to solve this problem, respectively. In the two methods, MIGA is used to search global optimum, and the sub-population size in MIGA is set as 200 for mathematical example. In interval optimization approach, the expansion order of chebyshev polynomials is set as 6. The sampling points and basis functions of mathematical example in each iteration of MIGA are shown in table 2, and the optimization results of two methods are shown in table 3 in which the validated values (bold in the table) of the objective is obtained by the scanning method in the uncertain range of the optimum point.

For the computation cost, the chebyshev tensor product (CTP) method is used to construct sixth-order chebyshev response surface, hence, 49 sampling points and 49 basis functions are required in each iteration of MIGA. In order to improve the approximation accuracy, the basis function in QSRS is increased to 90, and only 30 sampling points are required.

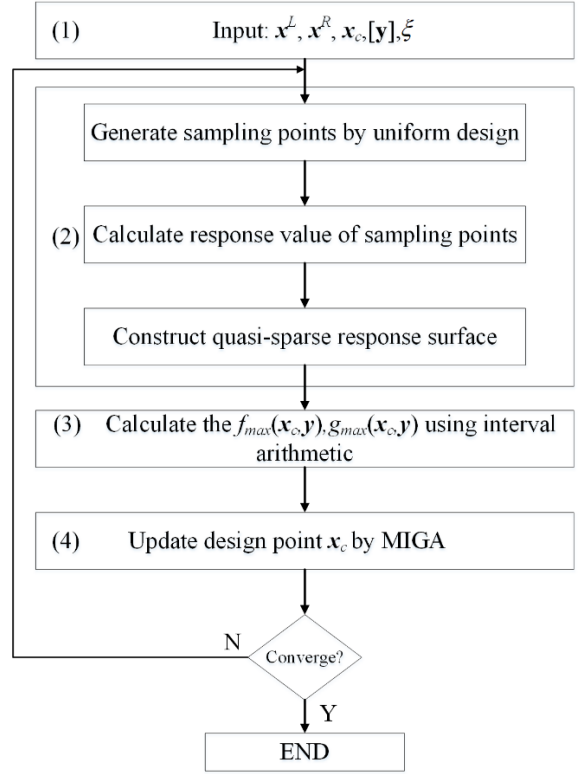

**FIGURE 4:** THE PROCESS OF QSRS APPROACH

Table 2. The sampling points and basis functions of mathematical example in each iteration of MIGA

|  | Interval optimization | QSRS approach |
|---|---|---|
| Sampling points | 49 | 30 |
| Basis functions | 49 | 90 |

Table 3. The optimization results of mathematical example

|  | Interval optimization | QSRS approach |
|---|---|---|
| $x_1$ | 1.7832 | 1.7901 |
| $x_2$ | -5.000 | -4.9933 |
| $f(x_1, x_2)$ | -4.3603(**-4.3603**) | -4.3332(**-4.3332**) |

In table 3, the objective $f(x_1, x_2)$ obtained by interval optimization approach and SRS approach equal to the validated value, which means that the interval arithmetic is effective in optimal design point obtained by the two method.

However, the interval optimization approach gives a local optimal solution which is -4.3603, because of approximation accuracy and wrapping effect in some design points updated in the iteration of MIGA. In contrary, the QSRS approach gives a better solution which is -4.3332.

### 4.2 Engineering example

A cylindrical vessel is capped at both ends by hemispherical heads as shown in Figure 5 [11]. $T_s(x_1)$ denotes thickness of the shell, $T_h(x_2)$ denotes thickness of the head, $R(x_3)$ denotes the inner radius, $L(x_4)$ denotes length of the cylindrical section of the vessel, not including the head. The objective of this problem is to minimize the total cost of material, forming and welding. Both four design variables are considered as interval design variables, and the width of them is 0.1.

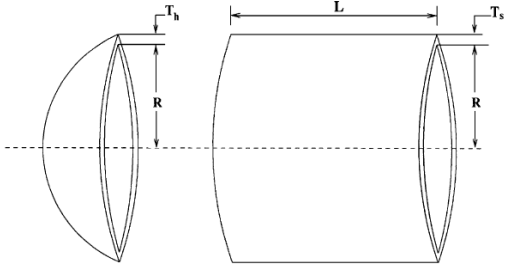

**FIGURE 5:** CENTER AND END SECTION OF PRESSURE VESSEL

$$\min_{x_{ic}} \max_{x_i \in [x_i]} f(x)$$
$s.t.\ g_1(x) = -x_1 + 0.0193x_3 \leq 0,$
$\quad g_2(x) = -x_2 + 0.00954x_3 \leq 0,$
$\quad g_3(x) = -\pi x_3^2 x_4 - \frac{4}{3}\pi x_3^3 + 1296000 \leq 0,$
$\quad g_4(x) = x_4 - 240 \leq 0,$
$\quad [x_i] = x_{ic} + [-\xi_i, \xi_i], \xi_i = 0.1, i = 1,2,3,4$
$\quad 1 + \xi_i \leq x_1, x_2 \leq 99 - \xi_i,$
$\quad 10 + \xi_i \leq x_3, x_4 \leq 200 - \xi_i.$

The interval optimization approach and QSRS approach are also used to solve this problem, respectively. The sub-population size in MIGA is set as 300 and the expansion order of chebyshev polynomials is set as 6 in interval optimization approach. The sampling points and basis functions of engineering example in each iteration of MIGA are shown in table 4, and the optimization results of two methods are shown in table 5 in which the validated values (bold in the table) of the objective is obtained by the scanning method in the uncertain range of the optimum point.

For the computation cost, chebyshev collocation method (CCM)) is employed to decrease the sampling points of sixth-order chebyshev response surface. Hence, 420 sampling points and 210 basis functions are required in each iteration of MIGA. However, only 100 sampling points is required in QSRS, as its basis function are increased to 300 to ensure the accuracy.

The optimization results in table 5 show that the interval optimization approach and QSRS approach give the same optimal point and all constrains are satisfied. However, the required sampling points of QSRS approach are reduced than interval optimization approach. It can be seen that the QSRS approach is effective for uncertainty structure design optimization, and as the dimension of the problem increases, the advantages of QSRS approach become more obvious.

Table 4. The optimization results of mathematical example

|  | Interval optimization | QSRS approach |
|---|---|---|
| Sampling points | 420 | 100 |
| Basis functions | 210 | 300 |

Table 5. The optimization results of mathematical example

|  | Interval optimization | QSRS approach |
|---|---|---|
| $x_1$ | 2.2254 | 2.2254 |
| $x_2$ | 1.2458 | 1.2458 |
| $x_3$ | 93.4970 | 93.4970 |
| $x_4$ | 100.2317 | 100.2317 |
| $g_1$ | -0.3190 | -0.3190 |
| $g_2$ | -0.2529 | -0.2529 |
| $g_3$ | -4.8606e6 | -4.8606e6 |
| $g_4$ | -139.6683 | -139.6683 |
| $f(x)$ | 4.6315e4(**4.6315e4**) | 4.6315e4(**4.6315e4**) |

## 5. CONCLUSION

In this paper, an efficient interval uncertainty optimization approach for uncertainty structure design optimization is proposed. Such approach combines QSRS with interval arithmetic to directly calculate the maximum value of objective, therefore, the inner loop optimization process is eliminated. With $l_1$ norm and $l_2$ norm penalty method and a few sampling points, the QSRS can be accurately constructed by a large number of basis functions. The results in mathematical example and engineering example show that the proposed approach gives a better performance than the interval optimization approach, meanwhile, less sampling points (30\100 for QSRS approach, 49\420 for interval optimization approach) in each iteration is required. For high-dimensional problems, the sampling points required in interval optimization approach will be very large, and the QSRS approach may be a better choice at this time.

## ACKNOWLEDGEMENTS
This research has been supported by National Natural Science Foundation of China under Grant Nos. 51775116, and 51374987, and 51405177, and NSAF U1430124.